\documentclass[a_4paper]{scrartcl}
\usepackage[OT2,T1]{fontenc}
\usepackage[utf8]{inputenc}
\usepackage{lmodern}
\usepackage{geometry}
\usepackage{amsmath}
\usepackage{amsfonts}
\usepackage{amssymb}
\usepackage{amsthm}
\usepackage{multirow}
\usepackage{hyperref}
\usepackage[all,knot]{xy}
\usepackage{rotating}
\usepackage{tabularx,caption}
\usepackage{bbm}
\xyoption{arc}

  \theoremstyle{plain}
 \newtheorem{thm}{Theorem}[section]

 \theoremstyle{definition}
 \newtheorem{defn}[thm]{Definition}
 \theoremstyle{definition}
 \newtheorem{defnthm}[thm]{Definition and Theorem}
 \theoremstyle{definition}
 \newtheorem{prop}[thm]{Proposition}
 \theoremstyle{definition}
 
 \theoremstyle{definition}
 
 \newtheorem{rem}[thm]{Remark}
 \theoremstyle{definition}
 
 \theoremstyle{definition}
 
 \theoremstyle{definition}
 
 \theoremstyle{definition}
 
 \theoremstyle{definition}
 
 \theoremstyle{definition}
 
  \theoremstyle{definition}

\DeclareMathOperator{\Gal}{Gal}
\DeclareMathOperator{\Aut}{Aut}

\usepackage[usenames,dvipsnames]{color}

\author{Ana Ros Camacho and Rachel Newton}

\title{Orbifold autoequivalent exceptional unimodal singularities}
\makeindex

\begin{document}
\maketitle

\abstractname{.~ We prove that different expressions of the same exceptional unimodal singularity are orbifold equivalent in the sense of \cite{CR,CRCR,RCN}. As in our previous paper \cite{RCN}, the matrix factorizations proving these orbifold equivalences depend again on certain parameters satisfying some equations whose solutions are permuted by Galois groups.}

\tableofcontents

\section{Introduction}

Thanks to the classification performed by Arnold in the late 60's, we know that unimodal singularities fall into 3 families of parabolic singularities, a three-suffix series of hyperbolic ones, and 14 exceptional cases (more details can be found in e.g. \cite{Ar,AGV}). A way to describe them is via regular weight systems \cite{Sai}, that is, a quadruple of positive integers $\left( a_1, a_2,a_3; h \right)$ with $a_1, a_2, a_3 < h$ ($h$ is called the \textit{Coxeter number}) satisfying:
\begin{itemize}
\item[--] $\mathrm{gcd} \left( a_1,a_2,a_3 \right)=1$, and
\item[--] there exists a \textit{potential} (meaning a polynomial $W \in \mathbb{C} \left[ x_1,x_2,x_3 \right]$  that has an isolated singularity at the origin) with the degrees of the variables $x_i$ being $\vert x_i \vert=\frac{2 a_i}{h} \in \mathbb{Q}_{\geq 0}$, $i \in \lbrace 1,2,3 \rbrace$ which is invariant under the Euler field $E$, that is, $$E.W=\left( \frac{a_1}{h} x_1 \frac{\partial}{\partial x_1}+\frac{a_2}{h} x_2 \frac{\partial}{\partial x_2}+\frac{a_3}{h} x_3 \frac{\partial}{\partial x_3} \right) W=W.$$ 
\end{itemize}

We say that a potential is \textit{homogeneous of degree $d\in\mathbb{Q}_{\geq 0}$} if in addition it satisfies  $$W \left( \lambda^{\vert x_1 \vert} x_1,\lambda^{\vert x_2 \vert} x_2,\lambda^{\vert x_3 \vert} x_3 \right)=\lambda^d W \left( x_1,x_2,x_3 \right)$$ for all $\lambda \in \mathbb{C}^{\times}$.  Requiring the potential $W$ to be invariant under the Euler vector field turns out to be equivalent to requiring $W$ to be homogeneous of degree 2. This argument goes as follows: a potential in three variables can only have seven possible shapes, as specified in \cite[Chapter 13]{AGV}, for example. Imposing invariance under the Euler field boils down to some conditions on the powers of the monomials in the potential. With the assignment of degrees made, one easily finds that these conditions are precisely the same as those we should impose if we want homogeneity of degree 2.

From now on, when we write `potential' we mean a potential that is homogeneous of degree 2. Potentials are the central concept of the present paper.

Throughout this paper, we will work over the graded ring $S:=\mathbb{C} \left[ x_1,x_2,x_3 \right]$ (although most of the definitions and results can be easily generalized to polynomial rings in $n$ variables over a field $\mathbbm{k} \subset \mathbb{C}$). 
We denote the set of all possible potentials with complex coefficients, and three variables, by $\mathcal{P}_\mathbb{C}$. To a potential $W \in \mathcal{P}_\mathbb{C}$, we can associate a number called the \textit{central charge}, which is defined as: 

\begin{equation}
c_W=\sum\limits_{i=1}^3 \left( 1-\vert x_i \vert \right)
\nonumber
\end{equation}

The central charge is related to the Coxeter number by the formula $c_W=\frac{h+2}{h}$ \cite{Mar}.\footnote{This formula is a special case of the following formula for potentials in $n$ variables: $$c_W=\frac{\left( n-2 \right) h-2 \epsilon_W}{h}$$ where $\epsilon_W:=\sum\limits_{i=1}^n a_i-h$ is the \textit{Gorenstein parameter} of $W$ \cite{Sai}. In order to recover the formula of \cite{Mar}, set $n=3$ and use the fact that exceptional unimodal singularities have Gorenstein parameter $\epsilon_W=-1$.}

The potential $W$ associated to a singularity is not necessarily unique. For the case of exceptional unimodal singularities, one can compute that $\lbrace E_{14},Q_{12},U_{12},W_{12},W_{13},Z_{13} \rbrace$ have multiple associated potentials, see Table \ref{Excepts1}. 

\begin{table}
\begin{center}
\begin{tabular}{c|c|c|c|c|c}
Type & Potential (v1) & Potential (v2) & Potential (v3) &  $\left( a_1, a_2,a_3; h \right)$ \\ \hline
$Q_{10}$ & $x^4+y^3+x z^2$ & --  & -- & $\left( 9,8,6;24 \right)$ \\
$Q_{11}$ & $x^3 y+y^3+x z^2$ & --  & -- & $\left( 7,6,4;18 \right)$ \\
$Q_{12}$ & $x^3 z+y^3+x z^2$ & $x^5+y^3+x z^2$ & --  & $\left( 6,5,3;15 \right)$ \\
$S_{11}$ & $x^4+y^2 z+x z^2$ &  --& -- & $\left( 5,4,6;16 \right)$ \\
$S_{12}$ & $x^3 y+y^2 z+x z^2$ & -- & -- & $\left( 4,3,5;13 \right)$ \\
$U_{12}$ & $x^4+y^3+z^3$ & $x^4+y^3+z^2 y$& $x^4+y^2 z+z^2 y$ & $\left( 4,4,3;12 \right)$ \\
$Z_{11}$ & $x^5+x y^3+z^2$ &  -- & -- & $\left( 8,6,15;30 \right)$ \\
$Z_{12}$ & $y x^4+x y^3+z^2$ &  --  & -- & $\left( 6,4,11;22 \right)$ \\
$Z_{13}$ & $x^3 z+x y^3+z^2$ & $x^6+y^3 x+z^2$  & -- & $\left( 5,3,9;18 \right)$ \\
$W_{12}$ & $x^5+y^2 z+z^2$ & $x^5+y^4+z^2$  & -- & $\left( 5,4,10;20 \right)$ \\
$W_{13}$ & $y x^4+y^2 z+z^2$ & $x^4 y+y^4+z^2$ & --  & $\left( 4,3,8;16 \right)$ \\
$E_{12}$ & $x^7+y^3+z^2$ & -- & -- & $\left( 14,6,21;42 \right)$ \\
$E_{13}$ & $y^3+y x^5+z^2$ &  --  & -- & $\left( 10,4,15;30 \right)$ \\
$E_{14}$ & $x^4 z+y^3+z^2$ & $x^8+y^3+z^2$ &  -- & $\left( 8,3,12;24 \right)$ \\ \hline
\end{tabular}
\caption{Exceptional unimodal singularities with associated potentials and regular weight systems. We use the labels v1, v2 and v3 to denote the different potentials associated to each singularity.} 
\label{Excepts1}
\end{center} 
\end{table}

In previous work (\cite{CRCR}), it was proven that potentials associated to simple singularities could be related by an equivalence relation which will be defined in the next subsection. In \cite{RCN}, we proved that potentials associated to $Q_{10}$ and $E_{14}$ -- strangely dual exceptional unimodal singularities -- could be related by this same equivalence relation. In this paper, we prove that the potentials associated to an exceptional unimodal singularity are equivalent in this sense.

\subsection*{On matrix factorizations}

One can define an equivalence relation between potentials as follows.

\begin{defn}
\begin{itemize}
\item Given a potential $W \in S$, a \textit{matrix factorization} of $W$ consists of a pair $\left( M, d^M \right)$ where
\begin{itemize}
\item $M$ is a $\mathbb{Z}_2$-graded free module over $S$;
\item $d^M \colon M \rightarrow M$ is a degree 1 $S$--linear endomorphism (the \textit{twisted differential}) such that:
\begin{equation}
d^M \circ d^M=W.\mathrm{id}_M.
\nonumber
\end{equation}
\end{itemize}
We may display the $\mathbb{Z}_2$-grading explicitly as $M=M_0 \oplus M_1$ and $d^M=\left( \begin{matrix} 0 & d_1^M \\ d_0^ M & 0 \end{matrix} \right)$.
If there is no risk of confusion, we will denote $\left( M,d^M \right)$ simply by $M$.

\item We call $M$ a \textit{graded matrix factorization} if, in addition, $M_0$ and $M_1$ are $\mathbb{Q}$-graded,  acting with $x_i$ is an endomorphism of degree $\vert x_i \vert$ with respect to the $\mathbb{Q}$-grading on $M$, and the twisted differential has degree 1 with respect to the $\mathbb{Q}$--grading on $M$\footnote{Note that these conditions imply that $W$ is homogeneous of degree 2 (as desired).}.
\end{itemize}
\end{defn}

We will denote by $\mathrm{hmf}^{\mathrm{gr}} \left( W \right)$ the idempotent complete full subcategory of graded finite--rank matrix factorizations: its objects are homotopy equivalent to direct summands of finite--rank matrix factorizations. The morphisms are homogeneous even (with respect to the $\mathbb{Z}_2$ degree) linear maps up to homotopy with respect to the twisted differential. This category is monoidal and has duals and adjunctions which can be described in a very explicit way. This leads to the following result which gives precise formulas for the left and right quantum dimensions of a matrix factorization.

\begin{prop}{\cite{CM,rigidity}}
Let $V \left( x_1,\ldots,x_m \right)$ and $W \left( y_1,\ldots,y_n \right)$ be two potentials and $M$ a matrix factorization of $W-V$. Then the left quantum dimension of $M$ is: $$\mathrm{qdim}_l \left( M \right)=\left( -1 \right)^{\binom{m+1}{2}} \mathrm{Res} \left[ \frac{\mathrm{str} \left( \partial_{x_1} d^M \ldots \partial_{x_m} d^M \partial_{y_1} d^M \ldots \partial_{y_n} d^M \right) d y_1 \ldots d y_n}{\partial_{y_1} W,\ldots,\partial_{y_n} W} \right]$$
and the right quantum dimension is:
$$\mathrm{qdim}_r \left( M \right)=\left( -1 \right)^{\binom{n+1}{2}} \mathrm{Res} \left[ \frac{\mathrm{str} \left( \partial_{x_1} d^M \ldots \partial_{x_m} d^M \partial_{y_1} d^M \ldots \partial_{y_n} d^M \right) d x_1 \ldots d x_m}{\partial_{x_1} V,\ldots,\partial_{x_m} V} \right].$$
\label{qdims}
\end{prop}

Quantum dimensions allow us to define the following equivalence relation:
\begin{defnthm}{\cite{CR,CRCR}}
Let $V$, $W$ and $M$ be as in the previous proposition. We say that $V$ and $W$ are \textit{orbifold equivalent} ($V \sim_{\mathrm{orb}} W$) if there exists a finite--rank matrix factorization of $V-W$ for which the left and the right quantum dimensions are non-zero. Orbifold equivalence is an equivalence relation in $\mathcal{P}_\mathbb{C}$.
\end{defnthm}

\begin{rem}{\cite[Proposition 6.4]{CR} (or \cite[Proposition 1.3]{CRCR}}) If two potentials $V$ and $W$ are orbifold equivalent, then their associated central charges are equal: $c_V=c_W$.
\label{charges}
\end{rem}

Let us give some comments on quantum dimensions and orbifold equivalences \cite{CRCR,CR}:
\begin{itemize}
\item \cite[Lemma 2.5]{CRCR} The quantum dimensions of graded matrix factorizations take values in $\mathbb{C}$. One can see this by counting degrees in the formulas given in Proposition \ref{qdims}.
\item The definitions of the quantum dimensions are also valid for ungraded matrix factorizations (in which case they will take values in $S$ instead of in $\mathbb{C}$). Furthermore, the quantum dimensions are independent of the $\mathbb{Q}$-grading on a graded matrix factorization.
\item So far, the difficulty of establishing an orbifold equivalence lies in constructing the explicit matrix factorization which proves it.
\end{itemize}

In \cite{RCN}, we proved that the potentials associated to the singularities $Q_{10}$ and $E_{14}$ (in either of its two expressions) are orbifold equivalent. Since orbifold equivalence is an equivalence relation, transitivity implies that the two potentials describing $E_{14}$ are also orbifold equivalent. The following question is hence a legitimate one: ``are different expressions of the same (exceptional unimodal) singularity orbifold equivalent to each other?'' \footnote{This question was originally posed by W. Ebeling, to whom we express our gratitude for stimulating this work.}

Potentials associated to the same singularity obviously share the same Coxeter number and central charge, which is a consequence of orbifold equivalence as noted in Remark \ref{charges}. Indeed we found the main result of this paper.

\begin{thm}
\label{conj}
We have the following orbifold equivalences between different potentials associated to the same exceptional unimodal singularities:
\begin{equation}
\begin{split}
Q_{12} \mathrm{(v1)} & \sim_{\mathrm{orb}} Q_{12} \mathrm{(v2)} \\
U_{12} \mathrm{(v1)} \sim_{\mathrm{orb}} & U_{12} \mathrm{(v2)} \sim_{\mathrm{orb}} U_{12} \mathrm{(v3)}
\\
W_{12} \mathrm{(v1)} & \sim_{\mathrm{orb}} W_{12} \mathrm{(v2)}
\\
W_{13} \mathrm{(v1)} & \sim_{\mathrm{orb}} W_{13} \mathrm{(v2)} \\
Z_{13} \mathrm{(v1)} & \sim_{\mathrm{orb}} Z_{13} \mathrm{(v2)} .
\end{split} 
\nonumber
\end{equation}
\end{thm}

Some further motivation for this work is as follows. In \cite[Section 1]{RCN}, we conjectured that potentials associated to strangely dual singularities (meaning different singularities that share the same Coxeter number) are orbifold equivalent. As mentioned above, we proved in \cite{RCN} that $E_{14}$ and $Q_{10}$ are indeed orbifold equivalent (independently of the potential we take to describe $E_{14}$), while orbifold equivalence for the pairs $\left( S_{11},W_{13} \right)$, $\left( Q_{11},Z_{13} \right)$ and $\left( Z_{11},E_{13} \right)$ remains to be proven. Some of these singularities have two possible potentials describing them, and proving that both descriptions are orbifold equivalent would allow us to use either of these potentials to pursue a proof of this conjecture. This may ease that task, which is work in progress \cite{NRC2}. In addition, for the interest of the first author, these equivalences may shed new insights into the Landau-Ginzburg/conformal field theory correspondence (LG/CFT) \cite{HW,LVW,VW,RC} for CFTs with central charge bigger than one.  It would be interesting to identify what should be the interpretation on the CFT side of these results and the role of Galois groups in them, {\`a} la \cite{Gep}. Due to computational difficulties, we postpone this analysis to \cite{NRC2}.

This paper is organized as follows. In the Section \ref{sec:pfs}, we describe the method obtained to generate the matrix factorizations proving the orbifold equivalences from Theorem \ref{conj},  give details for each of them, describe the equations the perturbation coefficients must satisfy and compute the Galois groups that permute them and any possible constraints from imposing non-vanishing of the quantum dimensions. For the sake of completeness, we also include a matrix factorization proving orbifold equivalence between the potentials describing $E_{14}$. Due to their size, the explicit descriptions of each matrix factorization are contained in Appendix \ref{appendix}. To finish, we sum up the results and outlook in Section \ref{Summary}.

\subsection*{Acknowledgements}
RN wishes to thank Bartosz Naskr\k{e}cki for a useful Mathematica tutorial. 
ARC's work is supported by the French-German research project \href{http://www.math.polytechnique.fr/SISYPH/sisyph-en.htm}{SISYPH} (programme blanc ANR-13-IS01-0001-01/02, DFG Program DFG No HE 2287/4-1, SE 1114/5-1). She also would like to warmly thank Andrew Taggart and Alex Pall for providing the soundtrack which supported the writing of this paper. 

\section{The proofs}
\label{sec:pfs}

All of the graded matrix factorizations (with $\vert x_i \vert=\frac{2 a_i}{h}$ given by the regular weight systems associated to each singularity) that we use to prove Theorem \ref{conj} are built in the same way, which we describe here. We start with a matrix factorization whose underlying module is $M=\mathbb{C} \left[ x,y,z,u,v,w \right]^{\oplus 8}$ and whose twisted differential $d^M=\left( d_{i j} \right)_{i,j \in \lbrace 1,\ldots,8 \rbrace}$ is constructed as follows. 
\begin{itemize}
\item[--] Take $d_{18}=d_{27}=d_{36}=d_{45}=d_{54}=d_{63}=d_{72}=d_{81}=0$.
\item[--] Impose \begin{equation}
\begin{split}
d_{64} &=d_{53}=-d_{28}=-d_{17} \\
d_{73} &=d_{62}=d_{48}=d_{15} \\
d_{74} &=-d_{52}=-d_{38}=d_{16} \\
d_{82} &=d_{71}=-d_{46}=-d_{35} \\
d_{83} &=-d_{61}=-d_{47}=d_{25} \\
d_{84} &=d_{51}=d_{37}= d_{26}
\end{split} 
\nonumber
\end{equation}\footnote{These constraints happened to be the same as those required in the matrix factorizations proving orbifold equivalence in \cite{RCN}.}
\item[--] All the remaining entries are zero.
\end{itemize} 

With this particular choice, we get a matrix factorization whose twisted differential squares directly to $-d_{17} d_{35} + d_{16} d_{25} - d_{15} d_{26}$. Hence, we will only need to specify the entries $\lbrace d_{15},d_{16},d_{17},d_{25},d_{26},d_{35} \rbrace$, which are listed in detail for each case in Appendix \ref{appendix}.

We take the $d_{15}$, $d_{16}$, $d_{17}$, $d_{25}$, $d_{26}$ and $d_{35}$ of the indecomposables $V_0$ in \cite{KST} associated to each singularity. Then we start perturbing each of these entries {\`a} la \cite{RCN,CRCR}, introducing in each entry monomials of the same total degree times some complex coefficient. Adjusting these coefficients so that the total twisted differential squares to the desired potential, we obtain the desired matrix factorization, together with several equations that these coefficients need to satisfy.

\subsection{$\mathbf{E_{14} \mathrm{(v1)} \sim_{\mathrm{orb}} E_{14} \mathrm{(v2)}}$}

Here we construct a matrix factorization of $u^8 + v^3 + w^2 -x^4 z - y^3 - z^2$ (i.e. $E_{14} \left( \mathrm{v2} \right)-E_{14} \left( \mathrm{v1} \right)$). A matrix factorization proving orbifold equivalence between these two potentials is specified in Equation \ref{E14} and it depends on the complex parameters: $$\lbrace a_1,a_2,a_3,a_4,b_1,b_2,b_3,c\rbrace.$$ The parameter $c$ must satisfy the following equation:
\begin{equation}
-1 - \frac{c^8}{4}=0
\label{conditionqdimE14}
\end{equation}
while the remaining 7 parameters are free.

We obtain the following left and right quantum dimensions:
\begin{itemize}
\item[] $\mathrm{qdim}_l \left( M \right)= -\frac{1}{8} \left( a_3 - b_3 + 4 c \right)$
\item[] $\mathrm{qdim}_r \left( M \right)= -\frac{c^7}{2}$
\end{itemize}


We need to impose on the values of $\lbrace a_3,b_3,c \rbrace$ some constraints in order to avoid zero values of the quantum dimensions. It is clear from Equation \ref{conditionqdimE14} that the right quantum dimension cannot vanish. To ensure non-vanishing of the left quantum dimension, we need to avoid the solutions of the system of equations formed by $a_3 - b_3 + 4 c=0$ and Equation \ref{conditionqdimE14}.

\subsubsection{Galois theory}
The equation $c^8+4=0$ factorises as $(c^4+2c^2+2)(c^4-2c^2+2)=0$. So we get two families of solutions for $c$, one for each irreducible factor. Each family consists of four solutions for $c$.

\smallskip

Family 1:  $c = \pm \sqrt{1\pm i}$

Family 2: $c = \pm \sqrt{-1\pm i}$

\smallskip

\noindent The field generated over $\mathbb{Q}$ by each family of solutions for $c$ is $\mathbb{Q}(\sqrt{1+ i},\sqrt{2})$ which has Galois group $D_8$ (dihedral of order $8$) over $\mathbb{Q}$. Family 1 and Family 2 are the two orbits for the action of the Galois group on the solutions for $c$.

\subsection{$\mathbf{Q_{12} \mathrm{(v1)} \sim_{\mathrm{orb}} Q_{12} \mathrm{(v2)}}$}

Here we construct a matrix factorization of $u^2 w + w^5 + v^3 - y^3 - x^3 z - x z^2$ (i.e. $Q_{12} \left( \mathrm{v2} \right)-Q_{12} \left( \mathrm{v1} \right)$). A matrix factorization proving orbifold equivalence between these two potentials is the one specified in Equation \ref{Q12} and it depends on seven complex parameters $\lbrace a_1,a_2,a_3,a_4,a_5,b_1, b_2 \rbrace$ which must satisfy the following six equations:
\begin{equation}
\begin{split}
b_1^2 + b_1 b_2 &=0 \\
-1+\left( -a_2 - a_4 + a_5 \right) b_1 b_2 - a_2 b_2^2 &=0 , \\
\left( a_2 a_3 + a_1 a_4 + a_2^2 a_4 - a_2 a_4^2 - a_1 a_5 + a_2 a_4 a_5 \right) \left( a_1 - a_3 - a_2 a_5 + a_4 a_5 - a_5^2 \right) &=0 , \\
b_1 \left( a_2 + a_4 - a_5 \right) \left( a_1 - a_3 - a_2 a_5 + a_4 a_5 - a_5^2 \right) &+ \\ b_2 \left[ a_2^2 \left( -a_4 - a_5 \right) + a_1 \left( a_2 - a_4 + a_5 \right) + a_2 \left(-2 a_3 + a_4^2 - a_5^2 \right)\right] &=0 , \\
\left( a_1 + a_2^2 + a_3 - a_4^2 + a_2 a_5 + a_4 a_5 \right) b_1 + \left( a_1 + a_2^2 - 2 a_2 a_4 + 2 a_2 a_5 \right) b_2 &=0 \\
a_1 a_3 + a_2^2 a_3 + 2 a_1 a_2 a_4 + a_2^3 a_4 - 2 a_2 a_3 a_4 - a_1 a_4^2 - 2 a_2^2 a_4^2 + a_2 a_4^3 - a_1 a_2 a_5 + 2 a_2 a_3 a_5 &+ \\ 
  + a_1 a_4 a_5+ a_2^2 a_4 a_5 - a_2 a_4^2 a_5 + a_2^2 a_5^2 - a_2 a_4 a_5^2 + a_2 a_5^3 &=0
\end{split}
\label{conditionsQ12} 
\end{equation}

We obtain the following left and right quantum dimensions:
\begin{equation}
\begin{split}
\mathrm{qdim}_l &=\frac{a_1}{10} \left(a_2 + a_4 - a_5\right) \left(2 a_3 + 3 a_5 \left(a_2 - a_4 + a_5\right) \right) b_1 \\ &+ \frac{a_1}{10}  \left(3 a_2^2 a_5 + a_2 \left( 2 a_3 + a_4 a_5 \right) + \left( 3 a_4 - 2 a_5 \right) \left( a_3 - a_4 a_5 + a_5^2 \right) \right) b_2\\ 
& - \frac{1}{10}\left( a_3 + a_5 \left( a_2 - a_4 + a_5 \right) \right)  \left( \left( a_2 + a_4 - a_5 \right) \left(2 a_3 + 3 a_5 \left( a_2 - a_4 + a_5 \right) \right) b_1 \right. \\ 
& \left.- \frac{a_2}{30} \left( a_3 + \left(a_4 - a_5 \right) \left(3 a_2 - 3 a_4 + 2 a_5\right) \right) b_2 \right) \\
\mathrm{qdim}_r &=-\frac{1}{5} \left( -a_4 b_1 + 5 a_5 b_1 + 3 a_5 b_2 + a_2 \left( b_1 + b_2 \right)\right)
\end{split}
\nonumber
\end{equation}

No solution of Equations \ref{conditionsQ12} makes the left quantum dimension vanish, but we need to impose on the values of $\lbrace a_1,a_2,a_3,a_4,a_5,b_1,b_2 \rbrace$ some constraints in order to avoid zero values of the right quantum dimension. The solutions of the Equations \ref{conditionsQ12} with $$-a_4 b_1 + 5 a_5 b_1 + 3 a_5 b_2 + a_2 \left( b_1 + b_2 \right)=0$$ are to be avoided.

\subsubsection{Galois theory}
Solving the equations using Mathematica shows that $a_5$ can be chosen freely and then the values for the other parameters lie in the field $\mathbb{Q}(a_5, \sqrt[5]{2}, \zeta_5, i)$ where $\zeta_5$ denotes a primitive $5$th root of unity. Generically, treating $a_5$ as an indeterminate, we get the Galois group 
\[\Gal(\mathbb{Q}(a_5, \sqrt[5]{2}, \zeta_5, i)/\mathbb{Q}(a_5))\cong \Gal(\mathbb{Q}(\sqrt[5]{2}, \zeta_5, i)/\mathbb{Q})\cong  C_2\times C_5\rtimes C_4\]
where the action of $C_4$ on $C_{5}$ is via an isomorphism $C_4\cong \Aut(C_{5})$.

\subsection{$\mathbf{U_{12} \mathrm{(v2)} \sim_{\mathrm{orb}} U_{12} \mathrm{(v3)}}$}
A matrix factorization of the potential $u^4+ v^3 + w^2 v- x^4 - y^2 z - y z^2$ (i.e. $U_{12} \mathrm{(v2)} -U_{12} \mathrm{(v3)}$) proving orbifold equivalence is specified in Equation \ref{U12v2v3}, and it depends on 6 complex parameters $\lbrace a_1,a_2,b_1,b_2,c_1,c_2 \rbrace$. The variables $\lbrace a_1,a_2,b_1,b_2 \rbrace$ have to satisfy the following equations:
 \begin{equation}
 \begin{split}
 -1 + a_1^2 b_1 - a_1 b_1^2 &=0 \\
 2 a_1 b_1 a_2 - b_1^2 a_2 + a_1^2 b_2 - 2 a_1 b_1 b_2 &=0 \\
 -1 + b_1 a_2^2 + 2 a_1 a_2 b_2 - 2 b_1 a_2 b_2 - a_1 b_2^2 &=0 \\
 a_2^2 b_2 - a_2 b_2^2 &=0
 \end{split}
 \label{condiUabdf}
 \end{equation}
while $c_1$ and $c_2$ can take any value. This matrix factorization has the following quantum dimensions:
\begin{equation}
\begin{split}
\mathrm{qdim}_r \left( M \right) &=-\frac{1}{24} \left( 8 c_1 a_2 - 4 a_1 b_2 + 8 c_1 b_2 + 4 b_1 \left( a_2 - 2 c_2 \right) + 8 a_1 c_2 \right) \\
\mathrm{qdim}_l \left( M \right)=& -\frac{1}{12} \left( 4  b_1^2 \left( -b_1 + c_1 \right) a_2 + 2 a_1^3 b_2 + a_1 b_1 \left( 2 c_1 \left( -a_2 + b_2\right) \right. \right. \\ 
&\left. \left. + b_1 \left( 3 a_2 + b_2 - 3 c_2\right)\right)+4 a_1^2 \left( c_1 b_2 + b_1 \left( 2 a_2 + 3 b_2 - 3 c_2\right) \right) \right)
\end{split}
\nonumber
\end{equation}

In order to avoid vanishing of the quantum dimensions, we have to restrict the values of $\lbrace a_1,b_1,c_1,a_2,b_2,c_2 \rbrace$ to only consider solutions of Equations \ref{condiUabdf} with
$$-8 c_1 a_2 - 4 a_1 b_2 + 8 c_1 b_2 + 4 b_1 \left( a_2 - 2 c_2 \right) + 8 a_1 c_2\neq 0$$ and 
\begin{equation}
\begin{split}
&-4  b_1^2 \left( -b_1 + c_1 \right) a_2 + 2 a_1^3 b_2 + a_1 b_1 \left( 2 c_1 \left( -a_2 + b_2\right) \right. \\ & \left. + b_1 \left( 3 a_2 + b_2 - 3 c_2\right)\right)+4 a_1^2 \left( c_1 b_2 + b_1 \left( 2 a_2 + 3 b_2 - 3 c_2\right) \right)\neq 0.
\end{split}
\nonumber
\end{equation}

\subsubsection{Galois theory}
All solutions for $a_1, a_2, b_1, b_2$ lie in the field $\mathbb{Q}(i, \zeta_3, \sqrt[3]{2})$, where $\zeta_3$ is a primitive $3$rd root of unity. The Galois group over $\mathbb{Q}$ is $S_3 \times C_2$. There are three families of solutions, corresponding to $a_2=0, b_2=0$ and $a_2=b_2$, respectively. Each family represents a distinct orbit for the action of the Galois group on the solutions for $a_1, a_2, b_1, b_2$.

\paragraph{\bf The family with $a_2=0$.} The equations give $b_1^3=1/2$, $a_1=2b_1$, $b_2^2=-b_1^2$. So $b_2=\pm ib_1$. There are $6$ solutions for $a_1,b_1,a_2,b_2$ in this family.

\paragraph{\bf The family with $b_2=0$.} The equations give $a_1^3=-1/2$, $b_1=2a_1$, $a_2^2=-a_1^2$. So $a_2=\pm ia_1$. There are $6$ solutions for $a_1,b_1,a_2,b_2$ in this family.

\paragraph{\bf The family with $a_2=b_2$. } The equations give $a_1^3=-1/2$, $b_1=-a_1$, $a_2^2=-a_1^2$. So $a_2=\pm ia_1$. There are $6$ solutions for $a_1,b_1,a_2,b_2$ in this family.

\subsection{$\mathbf{U_{12} \mathrm{(v1)} \sim_{\mathrm{orb}} U_{12} \mathrm{(v3)}}$}

A matrix factorization of $u^4+ v^3 + w^3 - x^4 - y^2 z - y z^2$ (i.e. $U_{12} \mathrm{(v1)}-U_{12} \mathrm{(v3)}$) proving orbifold equivalence is specified in
Equation \ref{U12v1v3} and depends on six complex parameters: $$\lbrace a_1,b_1,a_2,b_2,c_1,c_2 \rbrace$$

$\lbrace a_1,b_1,a_2,b_2 \rbrace$ must satisfy the equations:
 \begin{equation}
 \begin{split}
 -1 + a_1^2 b_1 - a_1 b_1^2 &=0 \\ 
 2 a_1 b_1 a_2 - b_1^2 a_2 + a_1^2 b_2 - 2 a_1 b_1 b_2 &= 0 \\
 b_1 a_2^2 + 2 a_1 a_2 b_2 - 2 b_1 a_2 b_2 - a_1 b_2^2 &= 0 \\
 -1 + a_2^2 b_2 - a_2 b_2^2 &= 0
 \label{conds1}
 \end{split}
 \end{equation}
$c_1$ and $c_2$ remain free. The quantum dimensions of this matrix factorizations are:
\begin{itemize}
\item[] $\mathrm{qdim}_l \left( M \right) = \frac{1}{9} \left( b_1 a_2 - 2 c_1 a_2 - a_1 b_2 + 2 c_1 b_2 + 2 a_1 c_2 - 2 b_1 c_2 \right)$
\item[] $\mathrm{qdim}_r \left( M \right) =4 (a_1^2 b_2 (4 a_2 - 3 b_2 + 2 c_2) +a_2 b_1 (3 a_2 b_1 - 2 b_1 b_2 - 2 a_2 c_1 + 4 b_2 c_1 - 2 b_1 c_2) +a_1 (-4 a_2^2 b_1 + 2 b_2 (b_1 b_2 + b_2 c_1 - 2 b_1 c_2) + a_2 (-4 b_2 c_1 + 4 b_1 c_2)))$
\end{itemize}

\smallskip

We are only interested in solutions to Equations \ref{conds1} for which both right and left quantum dimensions are nonzero.

\subsubsection{Galois theory}
 Solving the equations using Mathematica shows that all solutions for $a_1, b_1, a_2, b_2$ lie in $\mathbb{Q}(\zeta_3)$ for a primitive $3$rd root of unity $\zeta_3$. The Galois group of $\mathbb{Q}(\zeta_3)/\mathbb{Q}$ is $C_2$, generated by complex conjugation.

\subsection{$\mathbf{W_{12} \mathrm{(v1)} \sim_{\mathrm{orb}} W_{12} \mathrm{(v2)}}$}

A matrix factorization of $x^4+y^5+z^2-v^5 - u^2 w - w^2$ (i.e. $W_{12} \mathrm{(v2)}-W_{12} \mathrm{(v1)}$) proving the equivalence is given in Equation \ref{W12} and it depends on four complex parameters: $$\lbrace a_1,a_2,b_1,b_2 \rbrace$$

The parameters $\lbrace a_1,b_1,b_2 \rbrace$ must satisfy the following equations:
  \begin{equation}
  \begin{split}
& a_1^2 - a_1 b_1 =0 \\ 
& \frac{1}{4} \left(-4 - \left( 2 a_1 - b_1 + b_2 \right)^2 \left( \left( b_1 - b_2\right)^2 + 4 a_1 b_2 \right) \right)=0  
  \end{split}
  \label{conditionsW12}
  \end{equation}
  
  \quad This matrix factorization has quantum dimensions:
  \begin{equation}
  \begin{split}
  \mathrm{qdim}_r &= -\frac{1}{4} \left( 2 a_1 - b_1 + b_2 \right) \left( b_1^2 + 4 a_1 b_2 - 2 b_1 b_2 + b_2^2 \right)\\
   \mathrm{qdim}_l &=\frac{1}{2} \left( 2 a_1 - b_1 \right)
  \end{split}
  \nonumber
  \end{equation}
  
Non-vanishing of the right quantum dimension does not impose any additional constraints, since no solution of Equations \ref{conditionsW12} gives $\mathrm{qdim}_r\neq 0$. The solutions of Equations \ref{conditionsW12} with $\mathrm{qdim}_l= 0$ are as follows:
$a_1=b_1=0, b_2\in \lbrace \eta^{\pm 1} \sqrt{2},\eta^{\pm 3} \sqrt{2} \rbrace.$ These four solutions must be discarded.

\subsubsection{Galois theory}
There are two main families of solutions, as described below.

\paragraph{Family 1: $a_1=0$.} The second equation of \eqref{conditionsW12} can be written as $4+(b_1-b_2)^4=0$. Define a new variable $t=b_1-b_2$. Note that $(1+i)^4=-4$. Therefore, the solutions for $t$ come in two pairs: $1\pm i$ and $-1\pm i$. They lie inside $\mathbb{Q}(i)$, which has Galois group $C_2$ over $\mathbb{Q}$.

\paragraph{Family 2: $a_1=b_1.$} The second equation of \eqref{conditionsW12} simplifies to $4+(b_1+b_2)^4=0$. Let $u=b_1+b_2$. By the same argument as above, the solutions for $u$ are $1\pm i, -1\pm i$.

\paragraph{} So in both families, we can choose $b_2$ freely and then $a_1, b_1\in\mathbb{Q}(b_2, i)$.
Generically, treating $b_2$ as an indeterminate, we get the Galois group 
\[\Gal(\mathbb{Q}(b_2,i)/\mathbb{Q}(b_2))\cong \Gal(\mathbb{Q}(i)/\mathbb{Q})\cong C_2\]
generated by complex conjugation. 

\subsection{$\mathbf{W_{13} \mathrm{(v1)} \sim_{\mathrm{orb}} W_{13} \mathrm{(v2)}}$}

A matrix factorization of the potential $u^4 + v^4 u + w^2 - x^4 y - y^2 z - z^2$ (i.e. $W_{13} \mathrm{(v2)}$-$W_{13} \mathrm{(v1)}$) proving orbifold equivalence is described in Equation \ref{W13} and it depends on the complex parameters $$\lbrace  a_1,a_2,a_3,b,c,d,f,g \rbrace$$

We need to impose the following conditions on the parameters:
\begin{equation}
\begin{split}
-1-\frac{1}{4} \left( a_1 + b - c \right)^2 \left( 3 a_1 + 4 a_1 a_2 - b + 4 a_2 b + c - 4 a_2 c \right)^2 &=0 \\
-1 - \left( a_1 + b - c \right) d \left( a_3^2 - d + 2 a_3 f + f^2\right) &=0 \\
\left( a_1 + b - c \right) \left( a_3 + f \right) \left( a_3^2 - 2 d + 2 a_3 f + f^2 \right) &=0 \\
-2 \left( a_1 + b - c \right) \left( a_1 + 2 a_1 a_2 - b + 2 a_2 b + c - 2 a_2 c \right) &=0
\end{split}
\label{EqsW13}
\end{equation}

The quantum dimensions of this matrix factorization are:
\begin{equation}
\begin{split}
\mathrm{qdim}_l \left( M \right) &=-\frac{1}{8} \left( a_1 + b - c \right)^2 \left( a_1^2 \left( 4 + 6 a_2 \right) + \left( b - c \right) \left( \left( -1 + 4 a_2 \right) b + 2 c - 6 a_2 c \right) \right. \\ &\left. + a_1 \left( b + 10 a_2 b - 2 \left( c + 6 a_2 c \right) \right) \right) \left(a_3^3 + 3 a_3^2 f + f^3 - f g - a_3 \left( d - 3 f^2 + g \right) \right)
\end{split}
\nonumber
\end{equation}

$\mathrm{qdim}_r \left( M \right)=\frac{1}{16} \left( a_1 + 2 b - c \right) \left( 2 a_3 + f \right).$

\smallskip

We are only interested in solutions to Equations \ref{EqsW13} for which both right and left quantum dimensions are nonzero.

\subsubsection{Galois theory}
The equations above can be further simplified as follows:
\begin{itemize}
\item[] $-1 +d^2 \left( a_1 + b - c \right) =0 ,$
\item[] $ a_3 + f =0 ,$
\item[] $ a_1 + 2 a_1 a_2 - b + 2 a_2 b + c - 2 a_2 c =0$
\item[] $4+\left( a_1 + b - c \right)^4 =0$
\end{itemize}
Then the first equation gives $a_1 + b - c=d^{-2}$ and hence the last equation becomes $4+d^{-8}=0$. The third equation becomes 
$d^{-2}(2a_2-1)+2a_1=0$ so finally the equations reduce to

\begin{itemize}
\item[] $a_1 + b - c=d^{-2} ,$
\item[] $ a_3 + f =0 ,$
\item[] $d^{-2}(2a_2-1)+2a_1=0$
\item[] $4+d^{-8}=0 $.
\end{itemize}

The variables $a_2, a_3, b$ may be chosen freely. The others are determined by these and by $d$. Let $t=d^{-1}$ so $t^8+4=0$. The equation $t^8+4=0$ factorises as $(t^4+2t^2+2)(t^4-2t^2+2)=0$. So we get two families of solutions, one for each irreducible factor. Each family consists of four solutions for $t$.
\smallskip

Family 1:  $t = \pm \sqrt{1\pm i}$

Family 2: $t = \pm \sqrt{-1\pm i}$

\smallskip

The field generated over $\mathbb{Q}$ by each family of solutions for $t$ is $\mathbb{Q}(\sqrt{1+ i},\sqrt{2})$ which has Galois group $D_8$ (dihedral of order $8$) over $\mathbb{Q}$. Family 1 and Family 2 are the two orbits for the action of the Galois group on the solutions for $t$.

\subsection{$\mathbf{Z_{13} \mathrm{(v1)} \sim_{\mathrm{orb}} Z_{13} \mathrm{(v2)}}$}

A matrix factorization of the potential $u^6+v^3 u+w^2-x^3 z - z^2 - x y^3$ proving orbifold equivalence is described in Equation \ref{Z13}. It depends on the complex parameters $\lbrace a_1,a_2,a_3,b_1,b_2,c,d,f_1,f_2\rbrace$. The parameters $\lbrace c,d \rbrace$ must satisfy the following equations:
\begin{equation}
\begin{split}
-1 - \frac{c^6}{4} &=0\\
-1 + c d^3 &=0
\label{EqsZ13}
\end{split}
\end{equation}
The quantum dimensions of this matrix factorization are:
\begin{itemize}
\item[] $\mathrm{qdim}_l \left( M \right)=-\frac{c^6}{12} \left( 2 f_1 - 3 d f_2 \right)$
\item[] $\mathrm{qdim}_r \left( M \right)=-\frac{1}{9} \left( a_3 - b_1 + 3 c \right) \left( 3 d + f_2 \right).$
\end{itemize}

\smallskip

We are only interested in solutions to Equations \ref{EqsZ13} for which both right and left quantum dimensions are nonzero.

\subsubsection{Galois theory}
The equations give $d^3=1/c$ and $c^6=-4$. Let $t=d^{-1}$. Then $c=t^3$ (so $c$ is completely determined by $t$) and $t^{18}+4=0$. Let $t_0$ denote a solution to $t^{18}+4=0$. Then all solutions lie in the field $\mathbb{Q}(t_0, \zeta_9)$, where $\zeta_9$ denotes a primitive $9$th root of unity. 
We have 
\[\Gal(\mathbb{Q}(t_0, \zeta_9)/\mathbb{Q})=\Gal(\mathbb{Q}(t_0, \zeta_9)/\mathbb{Q}(\zeta_9))\rtimes \Gal(\mathbb{Q}(\zeta_9)/\mathbb{Q})\cong C_{18}\rtimes C_6\]

 Here $C_{18}$ is normal and $C_6$ acts on $C_{18}$ via an isomorphism $C_6\cong\Aut(C_{18})$. Explicitly, $\Gal(\mathbb{Q}(t_0, \zeta_9)/\mathbb{Q}(\zeta_9))$ is generated by $\sigma$, where $\sigma(t_0)=\zeta_{18}t_0=-\zeta_9t_0$ and $\sigma(\zeta_9)=\zeta_9$. And $\Gal(\mathbb{Q}(\zeta_9)/\mathbb{Q})$ is generated by $\tau$ where $\tau(\zeta_9)=\zeta_9^2$ and $\tau(t_0)=t_0$. We have $\tau\sigma\tau^{-1}=\sigma^{-7}$.

Note that $\Aut(C_{18})=\Aut(C_2\times C_9)=\Aut(C_2)\times \Aut(C_9)=\Aut(C_9)$. Therefore, 
\[C_{18}\rtimes \Aut(C_{18})=C_2\times C_9\rtimes \Aut(C_9)=C_2\times C_9\rtimes C_6.\]

\newpage
\section{Summary}
\label{Summary}

In Table \ref{SummaryTable}, we collect all the orbifold equivalences between potentials associated to exceptional unimodal singularities we have found to date, along with their respective Galois groups. We also include the number of equations we have, the number of variables in these equations, and the number of variables in the whole corresponding matrix factorization.

There are several points to remark here. As in \cite{RCN}, we keep observing a recurrent appearance of $C_2$ in the Galois groups. The Galois groups do not share the same order, and some are abelian while others are non-abelian. Some of them are repeated for different singularities, like $C_2$ (for $U_{12}$ and $W_{12}$) or $D_8$ (for $E_{14}$ and $W_{13}$). We lack a conceptual explanation for these apparent coincidences. Another limitation is the fact that the Galois groups do not take into account the extra constraints coming from the non-vanishing of the quantum dimensions.

Altogether, there are still many open questions concerning orbifold equivalences, and we hope these intermediate results may be early steps on a path towards a deeper understanding of this intriguing equivalence relation.

\begin{table}
\begin{center}
\begin{tabular}{c|c|c|c|c|} 
& \# Eqs & \# Vars & \# Vars & Galois \\ 
& & (in eq) & (in MF) & group \\ \hline
$E_{14} \mathrm{(v1)} \sim_{\mathrm{orb}} Q_{10}$ & 2 & 3 &  3 & $D_8 \times C_2$ \\
$E_{14} \mathrm{(v2)} \sim_{\mathrm{orb}} Q_{10}$ & 4 & 4 & 4 & $V_4$ \\
$E_{14} \mathrm{(v1)} \sim_{\mathrm{orb}} E_{14}\mathrm{(v2)}$ & 1 & 1 & 8 & $D_8$ \\
$Q_{12} \mathrm{(v1)} \sim_{\mathrm{orb}} Q_{12} \mathrm{(v2)}$ & 6 &7 & 7 & $C_2\times C_{5}\rtimes C_4$ \\
$U_{12} \mathrm{(v1)} \sim_{\mathrm{orb}} U_{12} \mathrm{(v2)}$ & 4 & 4 & 6 & $S_3 \times C_2$
\\
$U_{12} \mathrm{(v1)} \sim_{\mathrm{orb}} U_{12} \mathrm{(v3)}$ & 4 & 4 & 6 
 & $C_2$ \\
$W_{12} \mathrm{(v1)} \sim_{\mathrm{orb}} W_{12} \mathrm{(v2)}$ & 2 & 3 & 4 
& $C_2$ \\
$W_{13} \mathrm{(v1)}  \sim_{\mathrm{orb}} W_{13} \mathrm{(v2)}$ & 4 & 7 & 8 & $D_8$\\
$Z_{13} \mathrm{(v1)} \sim_{\mathrm{orb}} Z_{13} \mathrm{(v2)}$ & 2 & 2 & 9 & $C_2\times C_{9}\rtimes C_6$ \\ \hline 
\end{tabular}
\caption{Summary table. The action of $C_4$ on $C_{5}$ is via an isomorphism $C_4\cong \Aut(C_{5})$. The action of $C_6$ on $C_{9}$ is via an isomorphism $C_6\cong \Aut(C_{9})$. }
\label{SummaryTable}
\end{center} 
\end{table}

\quad

\newpage
\appendix
\section{Appendix}
\label{appendix}

Due to their size, we include in this appendix the explicit expressions of the matrix factorizations proving orbifold equivalence in each of the cases covered.

\subsection*{$\mathbf{E_{14}}$}
\begin{equation}
\begin{split}
d_{15} &=z+w + \kappa_2 u^4  +  a_1 x^4 + a_2 u^3 x+ a_3 u x^3 + a_4 u^2 x^2 ,\\
d_{16} &=v^2 + v y + y^2,\\
d_{17} &=x^3 z +\left( \frac{1}{2} a_2 \kappa_1 + b_2  \kappa_2  - c \left(a_2 b_2 - a_4 b_2 c - a_2 b_1 c + a_3 b_2 c^2 + a_2 b_3 c^2 + a_4 b_1 c^2 \right. \right. \\
& \left. \left. - a_2 a_1 c^3 - a_1 b_2 c^3 - a_4 b_3 c_1^3 - a_3 b_1 c^3 + a_4 a_1 c^4 + a_3 b_3 c^4 + a_1 b_1 c^4 - a_3 a_1 c^5 \right. \right. \\ 
&\left. \left. - a_1 b_3 c^5 + a_1^2 c^6 + a_4 \kappa_1- a_3 c\kappa_1+ \frac{1}{2} a_1 c^2 \left( 2\kappa_1 + b_1 \kappa_2 \right) - b_3 c \kappa_2 \right. \right. \\ 
&\left. \left. + a_1 c^2 \kappa_2 \right) \right) u^7 + \left( a_2 + b_2 - c \left( a_4 + b_1 - a_3 c - b_3 c + 2 a_1 c^2 \right) \right) u^3 w\\
&+\left( a_2 b_2 - a_4 b_2 c - a_2 b_1 c + a_3 b_2 c^2 + a_2 b_3 c^2 + a_4 b_1 c^2 - a_2 a_1 c^3 - a_1 b_2 c^3 \right. \\ & \left .- a_4 b_3 c^3 - a_3 b_1 c^3 + a_4 a_1 c^4 + a_3 b_3 c^4 + a_1 b_1 c^4 - a_3 a_1 c^5 - a_1 b_3 c^5 + a_1^2 c^6 + a_4 \kappa_1 \right. \\ 
& \left. -a_3 c \kappa_1+ a_1 c^2 \kappa_1 + b_1 \kappa_2 - b_3 c\kappa_2+ a_1 c^2 \kappa_2 \right) u^6 x+ \left( a_3 + b_3 - 2 a_1 c \right) u w x^2 \\ 
&+ \left( a_4 + b_1 - a_3 c - b_3 c + 2 a_1 c^2 \right) u^2 w x  + \left( a_3 a_1 + a_1 b_3 - a_1^2 c \right) u x^6 \\
&+ \left( a_4 b_2 + a_2 b_1 + a_3 \kappa_1 + b_3 \kappa_2 - c \left( a_3 b_2 + a_2 b_3 + a_4 b_1 - a_2 a_1 c - a_1 b_2 c \right. \right. \\ 
&- \left. \left. a_4 b_3 c - a_3 b_1 c + a_4 a_1 c^2 + a_3 b_3 c^2 + a_1 b_1 c^2 - a_3 a_1 c^3 - a_1 b_3 c^3 + a_1^2 c^4 
\right. \right. \\ 
& \left. \left. + a_1 \kappa_1 + a_1\kappa_2 \right) \right) u^5 x^2 + \left( -a_2 + b_2 - c \left( -a_4 + b_1 - \left( -a_3 + b_3 - c \right) c \right) \right) u^3 z \\ 
&+ \left( a_3 b_2 + a_2 b_3 + a_4 b_1 - a_2 a_1 c - a_1 b_2 c - a_4 b_3 c - a_3 b_1 c + a_4 a_1 c^2 + a_3 b_3 c^2 \right. \\ 
& \left. + a_1 b_1 c^2 - a_3 a_1 c^3 - a_1 b_3 c^3 + a_1^2 c^4 +a_1 \kappa_1+ a_1 \kappa_2 \right) u^4 x^3 + a_1^2 x^7+ 2 a_1 w x^3 \\ 
&+ \left( a_2 a_1 + a_1 b_2 + a_4 b_3 + a_3 b_1 - c \left( a_4 a_1 + a_3 b_3 + a_1 b_1 \right. \right. \\
& \left. \left. - c \left( a_3 a_1 + a_1 b_3 - a_1^2 c\right) \right) \right) u^3 x^4 + \left( -a_3 + b_3 - c \right) u x^2 z\\ 
& + \left( a_4 a_1 + a_3 b_3 + a_1 b_1 - c \left( a_3 a_1 + a_1 b_3 - a_1^2 c \right) \right) u^2 x^5  \\
& + \left( -a_4 + b_1 - \left( -a_3 + b_3 - c \right) c \right) u^2 x z  ,\\
d_{25} &=-v + y,\\
d_{26} &=-z+w+\kappa_1 u^4 + b_2 u^3 x + b_1 u^2 x^2 + b_3 u x^3 + a_1 x^4,\\
d_{35} &=x + c u,\\
\end{split}
\label{E14}
\end{equation}
with
\begin{equation}
\begin{split}
\kappa_1&=\frac{1}{2} \left( 2 b_2 c - 2 b_1 c^2 + 2 b_3 c^3 - c^4 - 2 a_1 c^4 \right)  \\ 
\kappa_2 &=a_2 c - b_2 c - a_4 c^2 + b_1 c^2 + a_3 c^3 - b_3 c^3 + c^4 + \kappa_1. 
\end{split}
\nonumber
\end{equation}

\subsection*{$\mathbf{Q_{12}}$}

\begin{equation}
\begin{split}
d_{15} &=z+ b_1 u + a_1 w^2 + a_2 w x ,\\
d_{16} &=v^2 + v y + y^2 ,\\
d_{17} &=x z+\left( a_2 b_1 + a_4 b_1 + a_2 b_2 \right) u w + \left( a_2 a_3 + a_1 a_4 + a_2^2 a_4 - a_2 a_4^2 + a_2 a_4 a_5 \right) w^3 \\ &+ a_2 a_4 w^2 x + a_5 w z ,\\
d_{25} &=-v+y ,\\
d_{26} &=-x z+a_5 b_2 u w + a_5 \left( -a_1 + a_3 + a_2 a_5 - a_4 a_5 + a_5^2 \right) w^3 + \left( b_1 + b_2 \right) u x \\ &+ a_3 w^2 x + a_4 w x^2 ,\\
d_{35} &=x^2+b_2 u + \left( -a_1 + a_3 + a_2 a_5 - a_4 a_5 + a_5^2 \right) w^2 + \left( -a_2 + a_4 - a_5 \right) w x .\\
\end{split}
\label{Q12}
\end{equation}

\subsection*{$\mathbf{U_{12}}$}
\subsubsection*{(v2) vs (v3)}
\begin{equation}
\begin{split}
d_{15} &=z+a_1 v + a_2 w , \\
d_{16} &=y z+ \left( -a_1^2 + a_1 b_1 + a_1 c_1 \right) v^2 + \left( -a_1 a_2 + b_1 a_2 + c_1 a_2 + a_1 \left( -a_2 + b_2 + c_2 \right) \right) v w \\ &+ a_2 \left(-a_2 + b_2 + c_2 \right) w^2 + c_1 v z + c_2 w z , \\
d_{17} &=u^3 + u^2 x + u x^2 + x^3 ,\\
d_{25} &=y+b_1 v + b_2 w ,\\
d_{26} &=y z+b_1 c_1 v^2 + \left( c_1 b_2 + b_1 c_2 \right) v w + b_2 c_2 w^2 + \left( -a_1 + b_1 + c_1 \right) v y + \left( -a_2 + b_2 + c_2 \right) w y ,\\
d_{35} &=x - u .\\
\end{split}
\label{U12v2v3}
\end{equation}

\subsubsection*{(v1) vs (v3)}

 \begin{equation}
 \begin{split}
 d_{15} &=z+ a_1 v + a_2 w ,\\
 d_{16} &= -y z+\left( -a_1^2 + a_1 b_1 + a_1 c_1 \right) v^2 + \left(-2 a_1 a_2 + b_1 a_2 + c_1 a_2 + a_1 b_2 + a_1 c_2 \right) v w \\ &+ a_2 \left( -a_2 + b_2 + c_2 \right) w^2 + c_1 v z + c_2 w z ,\\
 d_{17} &=  u^3 + u^2 x + u x^2 + x^3 ,\\
 d_{25} &= -y+ b_1 v + b_2 w ,\\
 d_{26} &= -y z+b_1 c_1 v^2 + \left( c_1 b_2 + b_1 c_2 \right) v w + b_2 c_2 w^2 - \left( -a_1 + b_1 + c_1 \right) v y \\ &- \left( -a_2 + b_2 + c_2 \right) w y ,\\
 d_{35} &= x - u . \\
 \end{split}
\label{U12v1v3}
 \end{equation}

\subsection*{$\mathbf{W_{12}}$}
\begin{equation}
\begin{split}
d_{15} &=z+w + a_1 u x + a_2 x^2 ,\\ 
d _{16} &=u w + b_1 x z + b_2 w x + \left( \left( -a_1 + b_1 \right) a_2 + a_1 \left( -a_2 + b_1 \left( -2 a_1 + b_1 - b_2 \right) \right) \right) x^3 ,\\
d_{17} &=v^4 + v^3 y + v^2 y^2 + v y^3 + y^4 ,\\ 
d_{25} &=u + \left( -2 a_1 + b_1 - b_2 \right) x ,\\
d_{26} &=z -w + \left( -a_1 + b_1 \right) u x + \left( -a_2 + b_1 \left( -2 a_1 + b_1 - b_2 \right) \right) x^2 ,\\
d_{35} &=-y+ v .\\
\end{split}  
\label{W12}
  \end{equation}

\subsection*{$\mathbf{W_{13}}$}

\begin{equation}
\begin{split}
d_{15} &=z+w+\frac{1}{2} \left( a_1^2 + 2 a_1^2 a_2 + 4 a_1 a_2 b - b^2 + 2 a_2 b^2 - 4 a_1 a_2 c + 
    2 b c - 4 a_2 b c - c^2 \right. \\ 
    & \left. + 2 a_2 c^2\right) u^2 + a_1 u y + a_2 y^2 , \\
d_{16} &= y z+\left(-\frac{3 a_1^3}{2} - 3 a_1^3 a_2 - a_1^3 a_2^2 - a_1^2 b - 6 a_1^2 a_2 b - 3 a_1^2 a_2^2 b + \frac{a_1 b^2}{2} - 3 a_1 a_2 b^2 - a_2^2 b^3  \right. \\ 
& \left. - 3 a_1 a_2^2 b^2  + \frac{3 a_1^2 c}{2} + 7 a_1^2 a_2 c + 3 a_1^2 a_2^2 c - a_1 b c + 8 a_1 a_2 b c + 6 a_1 a_2^2 b c - \frac{b^2 c}{2} + a_2 b^2 c \right. \\
&\left.  + 3 a_2^2 b^2 c + \frac{a_1 c^2}{2} - 5 a_1 a_2 c^2 - 
    3 a_1 a_2^2 c^2 + b c^2 - 2 a_2 b c^2 - 3 a_2^2 b c^2 - \frac{c^3}{2} + a_2 c^3 + a_2^2 c^3\right) u^3 \\ 
&- d g v^4 + \left(a_1 + 2 a_1 a_2 + 2 a_2 b + c - 2 a_2 c\right) u w + \left(-a_3 d - \left(-a_3 - f\right) g\right) v^3 x+ f v x^3  \\
& + \left(a_3^2 - d + a_3 f - g\right) v^2 x^2 + \left(a_1 c - \left(-a_1 - b + c\right) \left(a_1 a_2 + a_1 a_2^2 + a_2^2 b + a_2 c - a_2^2 c\right) \right.\\ 
&\left. + \frac{a_2}{2} \left(a_1^2 + 2 a_1^2 a_2 + 4 a_1 a_2 b - b^2 + 2 a_2 b^2 - 4 a_1 a_2 c + 2 b c - 4 a_2 b c - c^2 + 2 a_2 c^2\right) \right. \\
& \left. + a_2 \left(-a_1 \left(a_1 + 2 a_1 a_2 + 2 a_2 b + c - 2 a_2 c\right) - 
       b \left(a_1 + 2 a_1 a_2 + 2 a_2 b + c - 2 a_2 c\right)  \right. \right. \\ 
&\left. \left. + c \left(a_1 + 2 a_1 a_2 + 2 a_2 b + c - 2 a_2 c\right) + \frac{1}{2} \left(-a_1^2 - 2 a_1^2 a_2 - 4 a_1 a_2 b + b^2- 2 a_2 b^2 \right. \right. \right. \\ 
&\left. \left. \left.  + 4 a_1 a_2 c - 2 b c + 4 a_2 b c + c^2 - 2 a_2 c^2\right)\right)\right) u^2 y + 2 a_2 w y + a_2^2 y^3 + b u z \\
&  + \left(a_1 a_2 + a_1 a_2^2 + a_2^2 b + a_2 c - a_2^2 c\right) u y^2 , \\
 d_{17} &=  x^2 y+ \left(a_1 f + b f - c f\right) u v x + g v^2 y + a_3 v x y  +\left(a_1 a_3^2 + a_3^2 b - a_3^2 c - a_1 d - b d \right. \\ & \left. + c d + 2 a_1 a_3 f + 2 a_3 b f - 2 a_3 c f + a_1 f^2 + b f^2 - c f^2 - a_1 g - b g + c g\right) u v^2 ,\\
 d_{25} &=y+ \left( -a_1 - b + c \right) u , \\
 d_{26} &= -z+w+ c u y + a_2 y^2+\left(-a_1 \left(a_1 + 2 a_1 a_2 + 2 a_2 b + c - 2 a_2 c\right) -   b \left(a_1 + 2 a_1 a_2 \right. \right. \\
 &\left. \left.  + 2 a_2 b + c - 2 a_2 c\right) + c \left(a_1 + 2 a_1 a_2 + 2 a_2 b + c - 2 a_2 c\right) + \frac{1}{2} \left(-a_1^2 - 2 a_1^2 a_2 \right. \right. \\
 & \left. \left. - 4 a_1 a_2 b + b^2 - 2 a_2 b^2 + 4 a_1 a_2 c - 2 b c + 4 a_2 b c + c^2 - 2 a_2 c^2\right)\right) u^2 , \\
 d_{35} &= x^2+d v^2 + \left( -a_3 - f \right) v x .
\end{split}
\label{W13}
\end{equation}

\subsection*{$\mathbf{Z_{13}}$}

\begin{equation}
\begin{split}
d_{15} &= z+w+\left(a_2 c + b_2 c - a_3 c^2 - b_1 c^2 + 2 a_1 c^3 - \gamma \right) u^3 + a_1 x^3 + a_2 u^2 x + a_3 u x^2 ,\\
 d_{16} &= x y^2 + \left( -c f_1 + c d f_2 \right) u v^2 + \left( -f_1 - d \left( -d - f_2 \right) \right) v^2 x - c f_2 u v y + \left( -d - f_2 \right) v x y , \\
 d_{17} &= x^2 z + \left( a_2 \gamma + b_2 \left( a_2 c + b_2 c - a_3 c^2 - b_1 c^2 + 2 a_1 c^3  -\gamma \right)  \right. \\ 
 &\left. - c \left( a_2 b_2 + a_3 \gamma + b_1 \left( a_2 c + b_2 c - a_3 c^2 - b_1 c^2 + 2 a_1 c^3 - \gamma \right)\right. \right. \\ 
 & \left. \left. - c \left( a_3 b_2 + a_2 b_1 - a_2 a_1 c - a_1 b_2 c - a_3 b_1 c + a_3 a_1 c^2 + a_1 b_1 c^2 - a_1^2 c^3 +  a_1 \gamma  \right. \right. \right.\\ 
 &  \left. \left. \left. + a_1 \left( a_2 c + b_2 c - a_3 c^2 - b_1 c^2 + 2 a_1 c^3 -\gamma \right) \right) \right) \right) u^5 + \left( -d^3 + d f_1 - d^2 f_2 \right) v^3 \\ 
 & + \left( a_2 + b_2 - c \left( a_3 + b_1 - 2 a_1 c \right) \right) u^2 w + \left( a_3 + b_1 - 2 a_1 c \right) u w x + 2 a_1 w x^2 + a_1^2 x^5\\  
&+\left( a_2 b_2 +\gamma a_3 + b_1 \left( a_2 c + b_2 c - a_3 c^2 - b_1 c^2 + 2 a_1 c^3 -\gamma \right) \right. \\ 
&\left.  - c \left( a_3 b_2 + a_2 b_1 - a_2 a_1 c - a_1 b_2 c - a_3 b_1 c + a_3 a_1 c^2 + a_1 b_1 c^2 - a_1^2 c^3 + \gamma a_1 \right. \right. \\ 
& \left. \left.   + a_1 \left( a_2 c + b_2 c - a_3 c^2 - b_1 c^2 + 2 a_1 c^3 -\gamma \right) \right) \right) u^4 x\\ 
& + \left( a_3 b_2 + a_2 b_1 - a_2 a_1 c - a_1 b_2 c - a_3 b_1 c + a_3 a_1 c^2 + a_1 b_1 c^2 - a_1^2 c^3 + \gamma a_1 \right. \\ 
&\left. + a_1 \left( a_2 c + b_2 c - a_3 c^2 - b_1 c^2 + 2 a_1 c^3 + -\gamma_1\right) \right) u^3 x^2 \\
& + \left( a_2 a_1 + a_1 b_2 + a_3 b_1 - \left( a_3 a_1 + a_1 b_1 - a_1^2 c \right) \right) u^2 x^3 + \left( a_3 a_1 + a_1 b_1 - a_1^2 c\right) u x^4 \\ 
&  + f_1 v^2 y + f_2 v y^2 + \left( -a_2 + b_2 - \left( -a_3 + b_1 - c \right) c \right) u^2 z + \left( -a_3 + b_1 - c \right) u x z , \\
    d_{25} &=y + d v, \\
    d_{26} &=-z+w + a_1 x^3 + b_1 u x^2+ b_2 u^2 x+\gamma u^3  ,\\
    d_{35} &=x + c u ,
\end{split}
\label{Z13}
\end{equation}
where $\gamma=\frac{1}{2} \left( 2 b_2 c - 2 b_1 c^2 + c^3 + 2 a_1 c^3 \right)$.

\newpage

\newcommand\arxiv[2]      {\href{http://arXiv.org/abs/#1}{#2}}
\newcommand\doi[2]        {\href{http://dx.doi.org/#1}{#2}}
\newcommand\httpurl[2]    {\href{http://#1}{#2}}

\end{document}